\documentclass[a4paper,leqno]{article}
\usepackage{amssymb}
\usepackage{amscd}
\usepackage{amsmath}
\usepackage{amsthm}
\usepackage{mathrsfs}
\usepackage{graphics}
\usepackage{graphicx}
\usepackage{titlesec}
\titleformat{\section}[hang]
{\large\bf}
{\bf\thesection.}{.5em}{\bf }[ \vspace{+2pt}]
\titlespacing{\section}
{0pt}{1.5ex plus .1ex minus .2ex}{0pt}

\def\underset#1#2{{\mathrel{\mathop {{}_{} {#2}}\limits_{{#1}_{}}}}}
\def\upplim_#1{\underset{#1}{\overline\lim}\;}
\def\lowlim_#1{\underset{#1}{\underline\lim}\;}

\makeatletter
\def\@makefnmark{\hbox{\@textsuperscript{\normalfont
\@thefnmark)}}}
\renewenvironment{enumerate}
  {\ifnum \@enumdepth >3\relax\@toodeep\else
   \advance\@enumdepth\@ne
   \edef\@enumctr{enum\romannumeral\the\@enumdepth}%
   \list{\csname label\@enumctr\endcsname}{%
         \ifnum \@listdepth=\@ne \topsep.1\normalbaselineskip
           \else\topsep\z@\fi
         \parskip\z@ \itemsep\z@ \parsep\z@
         \labelwidth1pc \labelsep0.5pc
         \ifnum \@enumdepth=\@ne \leftmargin1pc\relax
           \else\leftmargin\leftskip\fi
         \advance\leftmargin +1.5pc
         \usecounter{\@enumctr}%
         \def\makelabel##1{\hss\llap{##1}}}%
   \fi}{\endlist}
\makeatother
\newtheorem{clm}[equation]{\indent \rm Claim}%

\newtheorem{prop}[equation]{Proposition}
\newtheorem{rmk}[equation]{\indent \rm {\it Remark}}

\newtheorem{thm}[equation]{Theorem}

\newtheorem{prob}[equation]{Problem}

\newcommand{\rB}{\mathrm{B}}
\newcommand{\C}{{\mathbf{C}}}

\newcommand{\co}{{\mathrm{co}}}

\newcommand{\fD}{{\mathfrak{D}}}

\newcommand{\del}{{\partial}}

\newcommand{\fa}{\forall}

\newcommand{\lto}{\longrightarrow}

\renewcommand{\O}{{\mathcal{O}}}

\newcommand{\PD}{{\mathrm{P}\hskip-1pt\Delta}}

\newcommand{\vp}{\varphi}

\newcommand{\R}{{\mathbf{R}}}

\numberwithin{equation}{section}

\title{
A Brief Proof of Bochner's Tube Theorem and\\ a Generalized Tube
}
\author{
J. Noguchi%
\thanks{
  Research supported in part by Grant-in-Aid
 for Scientific Research (C) 19K03511. 
}
}
\date{The University of Tokyo}
\begin{document}
\setlength{\baselineskip}{13pt}
\maketitle
\thispagestyle{empty}
\begin{abstract}
The aim of this note is firstly to give a new brief proof of classical
 Bochner's Tube Theorem (1938)
by making use of K. Oka's Boundary Distance Theorem (1942), showing
 directly that two points of
the envelope of holomorphy of a tube can be connected by a line
 segment. We then apply the same
idea to show that if an unramified domain
$\fD := A_1 + iA_2 \to \R^n + i\R^n = \C^n$
 with unramified real domains $A_j \to \R^n$
 is pseudoconvex, then both $A_j$ are univalent and
 convex (a generalization of Kajiwara's theorem).
 From the viewpoint of this result we discuss a
 generalization by M. Abe with
giving an example of a finite tube over $\C^n$ for which Abe's theorem no
 longer holds. The present
method may clarify the point where the (affine) convexity comes from.
\end{abstract}



\section{Introduction}
The following statement is classical and well-known as Bochner's Tube Theorem:
\begin{thm}[Bochner \cite{bo1}, \cite{bo2}, Stein \cite{st} ($n=2$)]
\label{bochtt}
Let $T_R=R+i\R^n$ be a tube (domain) of $\C^n$ with a domain (open, connected)
$R \subset \R^n$ as real base.
 Then the envelope of holomorphy of $T_R$
is $T_{\co(R)}$, where $\co(R)$ denotes the (affine)
convex hull of $R$.
\end{thm}

Our first aim is to give a new brief simple proof of this theorem,
based on Oka's Boundary Distance Theorem
(\S\ref{bochpf}). 

We then deal with a {\em generalized tube}
$\pi:A_1+iA_2 \to \C^n$ over $\C^n$ with real unramified
domains $\pi_j: A_j \to \R^n ~(j=1,2)$ and $\pi=\pi_1 +i \pi_2$.
 In \S\ref{gtpf} we also give another proof to Porten \cite{po}, Theorem 1.1:
\begin{thm}[Generalized tube]\label{gt}
If $A_1+iA_2$ is pseudoconvex, then the both
$A_j$ are univalent and convex subdomains of $\R^n$.
\end{thm}

The case where $A_j$ are univalent was obtained by
J. Kajiwara \cite{kaj63}, and the case where
$A_2=\R^n$ was dealt with by M. Abe \cite{abe}.
We will give counter-examples such that M. Abe's Theorem
does not holds for a {\em finite tube}; i.e., the part $A_2$
is  bounded in $\R^n$ (see  \S\ref{ftube}).

We will see the point where the (affine) convexity comes from
(see Remark \ref{pnt}).

\medskip
{\bf Acknowledgment.}
The author is very grateful to Professor Makoto Abe
for useful and helpful discussions during
the preparation of the present paper, and to Professors P. Pflug,
 P. Shapira and E. Porten
 for valuable informations on the present topics.

\section{Proof of Theorem \ref{bochtt}.}\label{bochpf} 
To be precise, a `domain' of $\R^n$ (or $\C^n$) is an open and connected subset
of $\R^n$ (or $\C^n$). If $X$ is a connected Hausdorff topological space
with a local homeomorphism $p: X \to \R^n$ (or $\C^n$), we call
$p: X \to \R^n$ (or $\C^n$) or simply $X$ ap
 {\em domain over} $\R^n$ (or $\C^n$).
If $p$ is injective, it is said to be univalent (schlicht)
 or otherwise multivalent
in general; a univalent domain over $\R^n$ (or $\C^n$) may be identified with
a domain of $\R^n$ (or $\C^n$).
In this paper, domains are always {\em unramified}.

For our proof we use the next two basic theorems:
 As for the envelope of holomorphy
we add the constructive existence for a convenience as an appendix
(cf.\ \S\ref{app} Appendix (1) at the end).

\begin{thm}\label{envh}
Every holomorphically separable domain $\fD$ over $\C^n$
admits an envelope of holomorphy, containing $\fD$ as a subdomain.
In particular, a univalent domain $\Omega$ of $\C^n$ admits
an envelope of holomorphy (multivalent in general),
 containing $\Omega$ as a subdomain.
\end{thm}

\begin{thm}[Boundary distance:
 Oka \cite{oka1}, \cite{oka2} VI (1942), IX (1953); \cite{hoer};
\cite{nog}]\label{oka}
If $\fD$ is a domain of holomorphy over $\C^n$,
then $-\log \delta (\zeta, \del\fD) ~(\zeta \in \fD)$
is a continuous plurisubharmonic function, where
$\delta (\zeta, \del\fD)$ denotes the distance function to the boundary
{\rm (cf.\ \S\ref{app} Appendix (2))}.
\end{thm}

Let $\pi: \hat T \to \C^n$ be the envelope of holomorphy of $T_R$
by Theorem \ref{envh}.
With $\hat R := \hat T \cap \pi^{-1} \R^n$, 
$\varpi=\pi|_{\hat R}:\hat R \to \R^n$ is a real domain
over $\R^n$ and $\varpi(\hat R) \subset \co(R)$.
Then $\hat T$ has a structure of a tube in the following sense:
\begin{equation}\label{tubetype}
\pi: \hat T= \hat R + i\R^n \lto \R^n +i\R^n= \C^n.
\end{equation}
It follows from Oka's Boundary Distance Theorem \ref{oka} that
$-\log \delta(\zeta, \del \hat T)$ is plurisubharmonic and satisfies
\begin{equation}\label{inv}
 -\log \delta(\zeta, \del \hat T)=
-\log \delta(\zeta+iy, \del \hat T), \quad
\fa y \in \R^n .
\end{equation}
With the local coordinates $\pi(p)=(x_j+iy_j)$, if $\delta(p, \del\hat T)$
 is of $C^2$-class, it satisfies the semi-positive definiteness:
\begin{equation}\label{hess}
\left(\frac{\del^2}{\del z_j \del \bar z_k}  -\log \delta(\zeta, \del
 \hat T)\right)_{j,k}
=\left(\frac{\del^2}{\del x_j \del  x_k}  -\log \delta(\zeta, \del
  \hat T)\right)_{j,k} \geq 0.
\end{equation}

We define a line segment $L[p,q] \subset \hat R$ connecting two
points $p, q \in \hat R$ as follows.
Let $L[\varpi(p), \varpi(q)] \subset \R^n$ be a line segment
connecting $\varpi(p)$ and $\varpi(q)$.
Then there is a unique connected component $L_p$ of the inverse
$\varpi^{-1}L[\varpi(p), \varpi(q)]$, containing $p$.
If $L_p \ni q$, we write $L_p=L[p,q] \subset \hat R$.
For mutually close $p,q \in \hat R$, $L[p,q]$ exists,
but in general  the existence is unknown at this moment.
If $p=q$, then $L[p,q]=\{p\}$ is considered as a special case
of  degenerate line segment.
Assuming the existence of $L[p,q]$, 
we see by \eqref{hess} that the restricted function
 $-\log \delta(\zeta, \del \hat T)|_{L[p,q]}$, even if it is not
 differentiable, is a convex
function on the line segment $L[p,q]$.
Therefore we have
\begin{equation}\label{max}
\min_{L[p,q]} \delta(\zeta, \del \hat T)=
\min_{\{p,q\}} \delta(\zeta, \del \hat T).
\end{equation}

\begin{clm} ~ If
 $S:=\left\{(p,q) \in {\hat R}^2 : 
 \exists\, L[p,q] \subset \hat R \right\} \subset {\hat R}^2 
$, then $S={\hat R}^2$. 
\end{clm}

Firstly, $S$ is non-empty and open.
It suffices to show that $S$ is closed in ${\hat R}^2$. 
Let $(p,q) \in \hat R^2$ 
 be an accumulation point of $S$.
Then there is a sequence of points
$(p_\nu,q_\nu) \in S ~(\nu=1,2, \ldots)$ such that
\[
 \lim_{\nu \to \infty} p_\nu=p, \quad
 \lim_{\nu \to \infty} q_\nu=q, \quad
L[p_\nu, q_\nu] \subset \hat R.
\]
By \eqref{max} there is a constant 
$\rho_0>0$ independent of $\nu$ such that the tubular neighborhood
$U_\nu$ (univalent) of every $L[p_\nu, q_\nu]$ with width $\rho_0$
is contained in $\hat R$.
Then for every sufficiently large $\nu$, $U_\nu \ni p,q$.
Therefore $L[p,q] \subset U_\nu \subset \hat R$; thus,
$(p,q) \in S$ and hence $S=\hat R^2$. 

It follows that $\varpi: \hat R \to \R^n$ is univalent.
For, otherwise, there were two points,
$p, q \in \hat R$ such that  $p \not=q$ and  $\varpi(p)=\varpi(q)$.
But there would be no line segment $L[p,q]$; contradiction.
Moreover, for arbitrary distinct  $p,q \in \hat R$,
$L[p,q] \subset \hat R$, and hence $\hat R$ is convex.
Thus, $\hat R=\co(R)$ and $\hat T=T_{\co(R)}$.
 \hfill $\square$\medskip

The above proof immediately implies the following generalization due to
 M. Abe \cite{abe}.
\begin{thm}\label{abth}
Let $\varpi: R \to \R^n$ be a real domain over $\R^n$
 and let $\pi: T_R=R+i\R^n \to \C^n$ be a
domain as in \eqref{tubetype}. Then, $T_R$ is a domain of holomorphy
if and only if $T_R$ is univalent and convex.
\end{thm}

\begin{rmk}\label{pnt}\rm
In the above proof, it was the point to deduce
the (affine) convexity  from \eqref{hess}, provided that the domain
is pseudoconvex or a domain of holomorphy.
\end{rmk}

\smallskip
{\bf  Notes.} Theorem \ref{bochtt} was proved by S. Bochner \cite{bo1}, \cite{bo2},
 and by K. Stein \cite{st} (Hilfssatz 1)
in $n=2$. Since then there have been many papers
 dealing with the proof (cf.\ Jarnicki--Pflug \cite{jp}, \S3.2
for more informations). The proofs were rather technically involved
(cf., e.g., \cite{bm} Chap.\ V, \cite{hoer} Chap.\ II).
The methods may be classified into four kinds:
\begin{enumerate}
\setlength{\itemsep}{-1.5pt}
\item
By Legendre polynomial expansions (Bochner \cite{bo2}, Bochner--Martin \cite{bm}).
\item
By a family of ellipses (Stein \cite{st} ($n=2$),
S. Hitotsumatsu \cite{hi},
 L. H\"ormander \cite{hoer} (Theorem 2.5.10), etc.)
\item
By the boundary distance function (H.J. Bremermann \cite{br}
in the case of  $n=2$).
\item
An approximation theorem of Bauendi-Treves (J. Hounie \cite{ho}).
\end{enumerate}

The present proof may belong to (iii) and was inspired by
Fritzsche--Grauert \cite{frgr} p.\,87 Exercise 1,
while in the textbook the notion of unramified domains is presented in the
subsequent section after it; so the supposed situation might be different
 to the present one.
It is also noticed that the observation of \eqref{hess} goes back
to Bremermann \cite{br} \S3.5.
 In the present proof as above, the univalence of
the envelope of holomorphy $\hat T$ and the convexity are proved at
the same time.

\section{Proof of Theorem \ref{gt}}\label{gtpf}
Let $\pi_j: A_j \to \R^n ~(j=1,2)$ be domains and
put
\[
\pi: \fD:=A_1+iA_2 \ni x+iy \to \pi_1(x)+ i\pi_2(y) \in \C^n.
\]
We call $\fD$ a {\em generalized tube (domain)}.
Let $y_0 \in A_2$ be arbitrarily fixed point, and take
a univalent ball neighborhood $\rB(y_0; 2\rho_0) \subset A_2$ with center $y_0$
and radius $2\rho_0>0$. The assumption implies that
the continuous function $\vp(z):=-\log\delta(z, \del\Omega)$ is
 plurisubharmonic in $\fD$.
Set
\[
 V=\{x \in A_1: \delta(x, \del A_1)<\rho_0 \}.
\]
Then the function $\vp(x+iy)$ in $x+iy \in V+i\rB_0(y_0; \rho_0)$
 is a function only in $x$. Therefore, $\vp(x+iy)=\vp(x+iy_0)$
is convex in $x\in V$. We set
\[
\psi(x)=\max\{\vp(x+iy_0), -\log \rho_0\}, \quad x \in A_1 .
\]
Then $\psi(x)$ is a continuous convex function in $A_1$.
The same arguments as in \S\ref{bochpf} with $\psi(x)$ imply that
$A_1$ is univalent and convex; the same is applied to $A_2$.
\hfill \qed

\section{Counter-examples of Abe's Theorem \ref{abth}
 for finite tubes}\label{ftube}
Here we give examples of `finite tubes' by replacing the imaginary part
$\R^n$ in Theorem \ref{abth} by a bounded domain, to say, an open ball,
for which the theorem no longer holds.

Let $0<R_1<R_2 \leq \infty$ and set
\begin{align*}
A &= \{x=(x_1, x_2) \in \R^2: R_1 <\|x\|:=(x_1^2+x_2^2)^{1/2} <R_2 \}, \\
\rB &= \{y=(y_1, y_2) \in \R^2 : \|y\| <R_1 \}.
\end{align*}
With complex coordinates $z_j=x_j+iy_j ~(j=1,2)$ we define
 a `finite tube' or a `tube of finite length' by
\[
 \Omega=A+i \rB \subset \C^2.
\]
We consider a holomorphic function $f(z)=z_1+iz_2 \in \O(\Omega)$
(it is the same with $f(z)=z_1 -iz_2$).
Since
\[
 |f(z)|=|x_1+ix_2 +i(y_1+iy_2)|\geq |x_1+ix_2|- |y_1+iy_2|>0,
\]
$g(z)=1/f(z) \in \O(\Omega)$; in particular, $g(z)$ is not holomorphic
at the origin $0$. Therefore we first note:
\begin{rmk}\rm
The envelope of holomorphy $\hat\Omega$ of $\Omega$ is {\em not}
 equal to $\co(A)+i \rB$.
 This give a counter-example for Kajiwara \cite{kaj63},
from which $\hat\Omega=\co(A)+i\rB$ should follow.
Cf.\ Jarnicki--Pflug \cite{jp}, \S3.3 for more examples and discussions.
\end{rmk}

 Let $ 2 \leq \nu \leq \infty$. For $ 2 \leq \nu < \infty$ we put
\begin{align*}
 A_\nu &= \left\{u=(u_1, u_2) \in \R^2: R_1^{1/\nu} < \|u\| <R_2^{1/\nu}\right\},\\
p_\nu &:  A_\nu \ni u=u_1+iu_2  \mapsto u^\nu=x_1+ix_2=(x_1,x_2)=x \in A,
\end{align*}
where the complex structures of `$u_1+iu_2$' and `$x_1+ix_2$'
are different and independent to that of $(z_1, z_2) \in \C^2$.
It follows that
$p_\nu$ is a local {\em real analytic} diffeomorphism between the annuli.
We put
\[
\pi_\nu:  \Omega_\nu=A_\nu \times \rB \ni (u, y) \to
p_\nu(u)+iy \in \Omega \hookrightarrow \C^2.
\]
Then $\pi_\nu:\Omega_\nu \to \C^2$ is a local real analytic diffeomorphism
and hence an unramified domain over $\C^2$.
We consider $f_\nu(z)=(f(z))^{1/\nu}=(x_1+ix_2+i(y_1+iy_2))^{1/\nu}$,
which is $\nu$-valued holomorphic in $z \in \Omega$. Note that
\[
 f_\nu(z)=(x_1+ix_2)^{1/\nu}\left(1+i \frac{y_1+iy_2}{x_1+ix_2}\right)^{1/\nu}:
\]
Here the latter product factor
 $\left(1+i \frac{y_1+iy_2}{x_1+ix_2}\right)^{1/\nu}$ has a
$1$-valued branch in $\Omega$, because 
\begin{equation}\label{les1}
\left| \frac{y_1+iy_2}{x_1+ix_2} \right|<1.
\end{equation}
Whereas the first factor $(x_1+ix_2)^{1/\nu}$ is defined to be $1$-valued
in $A_\nu$, and hence $f_\nu(z)$ is $1$-valued holomorphic in
 $\Omega_\nu$.
It follows that the domain $\pi_\nu: \Omega_\nu \to \C^2$ is holomorphically
 separable and $g_\nu=1/f_\nu \in \O(\Omega_\nu)$.

For $\nu=\infty$, we put
\[
\begin{array}{ccc}
 p_\infty: A_\infty=\{(u_1,u_2)\in \R^2: \log R_1 < u_1 < \log R_2,~ u_2
 \in \R\} & \to & A \\
\hbox{\rotatebox{90}{$\in$}} & \empty & \hbox{\rotatebox{90}{$\in$}}\\
 u=(u_1,u_2)  & \mapsto & e^{u_1}e^{iu_2}=(e^{u_1} \cos u_2,
 e^{u_1}\sin u_2) ~ .
\end{array}
\]
Then $p_\infty: A_\infty \to A$ is a local real analytic diffeomorphism.
Set
\[
\pi_\infty:  \Omega_\infty= A_\infty \times \rB \ni (u,y) \mapsto p_\infty(u)+iy \in
 \Omega \hookrightarrow \C^2.
\]
Then, $\pi_\infty:\Omega_\infty \to \C^2$ is an infinitely-sheeted
unramified domain over $\C^2$.

 We take $f_\infty(z)=\log f(z)$. Then we have
\[
 f_\infty (z)= \log (x_1+ix_2)+ \log\left(1+i
 \frac{y_1+iy_2}{x_1+ix_2}\right) ,
\quad z \in \Omega :
\]
Here, because of \eqref{les1} the second term
 $\log\left(1+i \frac{y_1+iy_2}{x_1+ix_2}\right)$
has a $1$-valued branch in $\Omega$ and the first term $\log(x_1+ix_2)$
is $1$-valued in $\Omega_\infty$, so that $f_\infty \in
\O(\Omega_\infty)$.
Therefore, the unramified domain $\pi_\infty:\Omega_\infty \to \C^2$
is holomorphically separable.
Since $f_\infty$ has no zero in $\Omega_\infty$,
$1/f_\infty \in \O(\Omega_\infty)$.

Thus we have:
\begin{prop}
Let the notation be as above.
For every $\nu$ with $2 \leq \nu \leq \infty$, $\pi_\nu: \Omega_\nu \to \C^2$
is a $\nu$-sheeted holomorphically separable unramified domain over $\C^2$, and
the envelope of holomorphy $\hat\pi_\nu :\hat\Omega_\nu \to \C^2$ of
 $\Omega_\nu$ is never univalent over  $\C^2$ and
 $\hat\pi_\nu(\hat\Omega_\nu) \not\ni 0$.
\end{prop}

We may propose at the end:
\begin{prob}\rm
Let $\Omega = A_1 +iB$ be a univalent generalized tube with
$A_1 \subset \R^n$ and an open ball $B \subset \R^n$.
\begin{enumerate}
\item
What is the envelope of holomorphy $\hat\Omega$ of $\Omega$?
\item
What is the condition of $A_1$ with which $\hat\Omega$ is univalent.
For example, if $A_1$ is simply connected or
contractible, is $\hat\Omega$ univalent?
\end{enumerate}
\end{prob}

\begin{rmk}\rm Very lately, Jarnicki--Pflug \cite{jp21}
dealt with the above problem for $\Omega=A_1+iB$ with
\[
 A_1=\{x \in \R^n : R_1 < \|x\| < 1\}, \quad
B=\{y \in \R^n: \|y\|<R_2  \}.
\]
This case is interesting in view of the above counter-example, and
they proved that $\hat\Omega$ is univalent and given by
\[
 \hat\Omega= \{x+iy \in \R^n+i\R^n: \|x\|<1, \|y\|<R_2,
\|y\|^2 <\|x\|^2+R_2^2-R_1^2.\}
\]
\end{rmk}

\section{Appendix}\label{app}
{\bf (1) Envelope of holomorphy.}
In quite a few references, the notion of the envelope of holomorphy
of domains over $\C^n$ are presented in a rather
sophisticated manner. For our aim the following
 simple-minded constructive
existence is sufficient.

We first fix a notation.
If $\fD$ is a connected Hausdorff space and $\pi: \fD \to \C^n$
is a local homeomorphism, $\pi: \fD \to \C^n$ or simply $\fD$ is
called a (unramified Riemann) {\em domain over} $\C^n$.
If $\pi$ is injective, $\fD$ is said to be univalent.
A domain $\fD$ over $\C^n$ naturally admits a structure of complex manifold
such that $\pi$ is a local biholomorphism;  the set of all holomorphic
functions on $\fD$ is denoted by $\O(\fD)$.

For an element $f \in \O(\fD)$ and  a point $p \in \fD$ 
there is a small polydisk neighborhood of $a=\pi(p)$
which is identified with a neighborhood of $p$, and
$f$ is written there as a convergent power series in the local coordinate
$z$:
\[
\underline{f}_{p}:=f(z) =\sum_{\alpha} c_\alpha (z-a)^\alpha.
\]
If for  two points $p,q \in \fD$ with $p\not=q$ and $\pi(p)=\pi(q)$
there is an element $f \in \O(\fD)$ such that
$\underline{f}_p \not= \underline{f}_q$, then
$\pi: \fD \to \C^n$ is said to be {\em holomorphically separable}.

We fix a point $p_0 \in \fD$. We consider a curve $C^b$ in $\C^n$
with the initial point $a=\pi(p_0)$ and the end point $b \in \C^n$
such that every analytic function $\underline{f}_{p_0}$ at $a$
defined by $f \in \O(\fD)$ can be analytically continued along $C^b$,
and defines an analytic function, 
denoted by $f_{C^b}(z)$, at the end point $b$.
Let $\Gamma$ denote the set of all such curves $C^b$.
If $C^b, C'^b \in \Gamma$ are homotope through
 a continuous family of curves belonging to $\Gamma$,
then $\underline{f_{C^b}}_b = \underline{f_{C'^b}}_b$.
We denote by $\{C^b\}$ the homotopy class in the above sense,
and write $\underline{f_{\{C^b\}}}_b :=\underline{f_{C^{b}}}_b$.

We fix a polydisk $\PD \subset \C^n$ with center at the origin.
For $f \in \O(\fD)$ and $C^b \in \Gamma$
there is a polydisk neighborhood
$b+ r\PD ~(r>0)$ of $b$ where
$\underline{f_{\{C^b\}}}_b (z)$ converges.
Let $r(\{C^b\}, f)$ be the supremum of such $r$, and
let  $\Gamma^\dagger$ denote all of $\{C^b\}$ such that
$\inf_{f \in \O(\fD)} r(\{C^b\}, f)>0$.

For two element $\{C^b\},\, \{C'^{b'}\}$ of $\Gamma^\dagger$
we define an equivalence relation 
$\{C^b\} \sim \{C'^{b'}\}$ by
\[
 b=b',  ~\quad \underline{f_{\{C^b\}}}_b = \underline{f_{\{C'^{b'}\}}}_{b'} ,
\quad \fa f \in \O(\fD) .
\]
Let $[\{C^b\}]$ stand for the equivalence class, and let
\[
 \hat\fD=\Gamma^\dagger /\sim,
\quad \hat\pi: [\{C^b\}] \in \hat\fD \to b \in \C^n
\]
be respectively the quotient set and the natural map.
It follows from the construction that
$\hat\pi: \hat\fD \to \C^n$ gives rise to a holomorphically separable
(unramified) domain over $\C^n$.
Since $\fD$ is arc-wise connected, $\hat\fD$ is independent
of the choice of $p_0 \in \fD$.
There is a natural holomorphic map
$\eta: \fD \to \hat\fD$ with $\pi=\hat\pi \circ \eta$.
If $\fD$ is holomorphically separable, then
$\eta$ is an inclusion map and
$\fD$ is a subdomain of $\hat\fD$.

We call $\hat\pi: \hat\fD \to \C^n$ the {\em envelope of holomorphy} of
$\fD$.
In the case of $n \geq 2$, even if $\fD$ is univalent,
the envelope of holomorphy $\hat\fD$ of $\fD$ may be
(infinitely) multi-sheeted over $\C^n$ in general.
If $\eta: \fD \to \hat\fD$ is biholomorphic
$(\fD=\hat\fD)$, $\fD$ is called a {\em domain of holomorphy}.

{\bf (2) Boundary distance.}
The boundary distance $\delta(\zeta, \del\fD)$ is defined as follows.
For a point $\zeta \in \fD$ there is an open ball
$B(\pi(\zeta);r)\subset\C^n$ with center $\pi(\zeta)$
and radius $r ~(>0)$ such that the connected component $U(\zeta; r)$
of $\pi^{-1}B(\pi(\zeta); r)$ containing $\zeta$ is
biholomorphically mapped onto $B(\pi(\zeta); r)$ by $\pi$.
We write $\delta(\zeta, \del \fD)$ for the supremum of such $r$,
which is called the {\em boundary distance}.
The proof of Theorem \ref{oka} is similar to the case of univalent
domains.

In place of an open ball we may use a polydisk $\PD$ with center
at $0$ in the above definition. Then the boundary distance is
denoted by $\delta_\PD(\zeta, \del\fD)$; Theorem \ref{oka} holds
with `$- \log_\PD(\zeta, \del\fD)$'.

{

}

\bigskip
\begin{flushright}
Junjiro Noguchi\\
Graduate School of Mathematical Sciences\\
University of Tokyo\\
Komaba 3-8-1, Meguro-ku, Tokyo 153-8914, Japan\\
E-mail address: noguchi@ms.u-tokyo.ac.jp
\end{flushright}

\end{document}